\documentclass[11pt]{amsart}
\usepackage{amsmath, amsthm, latexsym, amssymb, amsfonts, epsfig, epsf}

\newenvironment{pf}{\proof[\proofname]}{\endproof}
\theoremstyle{plain}
\newtheorem{Th}{Theorem}[section]

\newtheorem{Prop}[Th]{Proposition}
\newtheorem{Lemma}[Th]{Lemma}
\numberwithin{equation}{section} \theoremstyle{definition}
\newtheorem{Rem}[Th]{Remark}

\newcommand{\cal}[1]{\mathcal{#1}}

\newcommand{\mS}{\mathbb S}

\newcommand{\R}{\mathbb R}

\newcommand{\D}{\Delta}

\newcommand{\la}{\langle}
\newcommand{\ra}{\rangle}

\newcommand{\supp}{\operatorname{supp}}

\newcommand{\bd}{\partial}
\newcommand{\relint}{\operatorname{relint}}

\newcommand{\rl}[1]{Lemma~\ref{L:#1}}
\newcommand{\rp}[1]{Proposition~\ref{P:#1}}

\newcommand{\re}[1]{(\ref{e:#1})}

\newcommand{\rt}[1] {Theorem~\ref{T:#1}}

\begin{document}

\title[Characterization of simplices via the Bezout inequality]{Characterization of simplices via the Bezout inequality for mixed volumes}



\author{Christos Saroglou}
\address[Christos Saroglou]{Department of Mathematical Sciences\\ Kent State University\\ Kent, OH USA}
\email{csaroglo@math.kent.edu}
\author{Ivan Soprunov}
\address[Ivan Soprunov]{Department of Mathematics\\ Cleveland State University\\ Cleveland, OH USA}
\email{i.soprunov@csuohio.edu}
\author{Artem Zvavitch}
\address[Artem Zvavitch]{Department of Mathematical Sciences\\ Kent State University\\ Kent, OH USA}
\email{zvavitch@math.kent.edu}
\thanks{The third author is supported in part by U.S. National Science Foundation Grant DMS-1101636 and by the  Simons Foundation.}
\keywords{Convex Bodies, Mixed Volume, Convex Polytopes, Bezout Inequality, Aleksandrov--Fenchel Inequality}
\subjclass[2010]{Primary 52A39, 52A40, 52B11}

\date{}

\maketitle

\begin{abstract}
We consider the following Bezout inequality for mixed volumes:
$$V(K_1,\dots,K_r,\D[{n-r}])V_n(\D)^{r-1}\leq \prod_{i=1}^r V(K_i,\D[{n-1}])\ \text{ for }2\leq r\leq n.$$
It was shown previously that the inequality is true for any $n$-dimensional simplex $\D$ and any convex
bodies $K_1, \dots, K_r$ in $\R^n$. It was conjectured that simplices are the only convex bodies 
for which the inequality holds for arbitrary bodies $K_1, \dots, K_r$ in $\R^n$. In this paper we prove that
this is indeed the case if we assume that $\D$ is a convex polytope. Thus the Bezout inequality characterizes simplices in the class of convex $n$-polytopes.
In addition, we show that if a body $\D$ satisfies the Bezout inequality for all  bodies $K_1, \dots, K_r$ 
then the boundary of $\D$ cannot have {\it strict} points. In particular, it cannot have points with positive Gaussian curvature.
\end{abstract}


\section{Introduction}

It was noticed in \cite{SZ} that the classical Bezout inequality in algebraic geometry \cite[Sec. 8.4]{F} 
together with the  Bernstein--Kushnirenko--Khovanskii bound \cite{Be, Kush, Kho} produces a
new inequality involving mixed volumes of convex bodies:
\begin{equation}\label{e:BI0}
V(K_1,\dots,K_r,\D[{n-r}])V_n(\D)^{r-1}\leq \prod_{i=1}^r V(K_i,\D[{n-1}])\ \text{ for }2\leq r\leq n.
\end{equation}
Here $\Delta$ is an $n$-dimensional simplex and $K_1,\dots, K_r$ are arbitrary convex bodies in~$\R^n$.
Throughout the paper $V_n(K)$ denotes the $n$-dimensional Euclidean volume of a body $K$ and
$V(K_1,\dots, K_n)$ denotes the $n$-dimensional mixed volume of bodies $K_1,\dots, K_n$. Furthermore,
$K[{m}]$ indicates that the body $K$ is repeated $m$ times in the expression for the mixed volume.

In \cite{SZ} it was conjectured that the Bezout inequality characterizes simplices, that is if $\D$ is a convex body such
that \re{BI0} holds for all convex bodies $K_1,\dots, K_r$ then $\D$ is necessarily a simplex 
(see \cite[Conjecture 1.2]{SZ}). 
It was proved that $\D$ has to be indecomposable (see \cite[Theorem 3.3]{SZ}) which, in particular, confirms
the conjecture in dimension $n=2$. In the present paper we prove this conjecture for the class of 
convex polytopes.
\begin{Th}\label{T:simplex0} Fix $2\leq r\leq n$.
Let $\D$ be a convex $n$-dimensional polytope in $\R^n$ satisfying \re{BI0} for
all convex bodies $K_1,\dots, K_r$ in $\R^n$. Then $\D$ is a simplex.
\end{Th}
Although the above theorem covers a most natural class of convex bodies,  in full generality the conjecture remains
open. Going outside of the class of polytopes we show that if a convex body $\D$ satisfies \re{BI0} for
all convex bodies $K_1,\dots, K_r$ in $\R^n$ then $\D$ cannot have 
strict points. We say a boundary point $x\in K$ is a {\it strict point}
if $x$ does not belong to any segment contained in the boundary of $K$.

\begin{Th}\label{T:last0} Fix $2\leq r\leq n$.
Let $\D$ be an $n$-dimensional convex body  in $\R^n$ satisfying \re{BI0} for
all convex bodies $K_1,\dots, K_r$ in $\R^n$. Then $\D$ does not contain any strict points. 
\end{Th}
In particular, we see that $\D$ 
cannot have points with positive Gaussian curvature.

Let us say a few words about the idea behind the proofs of Theorems \ref{T:simplex0} and \ref{T:last0}. First,
note that it is enough to prove the theorems in the case of $r=2$ as this implies the general statement. Thus we are going to restate \re{BI0} for $r=2$ as
follows
\begin{equation}\label{e:BI-2-0}
V(L,M,K[n-2])V_n(K)\leq V(L,K[n-1])V(M,K[n-1]),
\end{equation}
where $L$ and $M$ are convex bodies and $K$ is a polytope.
The fact that there is equality in \re{BI-2-0} when $L=K$ allows us to see this as a variational problem, by fixing an appropriate body $M$ and using an appropriate deformation $L=K_t$ of $K$. In the case of \rt{simplex0}, $K_t$ is obtained from $K$ by moving one of its facets along the direction of its normal unit vector. In the case of \rt{last0}, $K_t$ is obtained from $K$ by cutting out a small cup in a neighborhood of a strict point.


\section{Preliminaries}
In this section we collect basic definitions and set up notation. As a general reference on the theory of 
convex sets and mixed volumes we use R. Schneider's book ``Convex bodies: the Brunn-Minkowski theory" \cite{Sch1}.

A {\it convex body} is a non-empty convex compact set. 
A {\it (convex) polytope} is the convex hull of a finite set of points. An $n$-dimensional polytope
is called an {\it $n$-polytope} for short.
For $x,y\in\R^n$ we write $\la x, y\ra$ for the inner product of $x$ and $y$. We use $\mS^{n-1}$ to
denote the $(n-1)$-dimensional unit sphere and $B(x,\delta)$ to denote the closed Euclidean ball of radius
$\delta>0$ centered at $x\in\R^n$.

For a convex body $K$ the function $h_K:\mS^{n-1}\to\R$, 
$h_K(u)=\max\{\la x, u\ra\ |\ x \in K\}$ is the {\it support function} of $K$.
For every $u\in\mS^{n-1}$  we write $H_K(u)$ to denote the
supporting hyperplane for $K$ with outer normal $u$
$$H_K(u)=\{x\in\R^n : \la x,u\ra = h_K(u)\}.$$ 
Furthermore, we use $K^u$ to denote the face $K\cap H_K(u)$ of $K$.

Let $\beta$ be a subset of the boundary $\bd K$ of a convex body $K$. The {\it spherical image} $\sigma(K,\beta)$ of 
$\beta$ with respect to $K$ is defined by
$$\sigma(K,\beta)=\{u\in \mS^{n-1}: \exists x\in\beta,\textnormal{ such that }\langle x,u\rangle=h_K(u)\}.$$
If $\Omega$ is a subset of $\mS^{n-1}$ define the {\it inverse spherical image} $\tau(K,\Omega)$ of $\Omega$ with respect to $K$ by
$$\tau(K,\Omega)=\big\{x\in\bd K:\exists u\in\Omega,\textnormal{ such that }\langle x,u\rangle=h_K(u)\big\}.$$

The {\it surface area measure} $S(K,\cdot)$ of $K$ (viewed as a measure on $\mS^{n-1}$) is defined as
$$S(K,\Omega)={\cal{H}}^{n-1}\big(\tau(K,\Omega)\big),\ \textnormal{ for }\Omega\textnormal{ a Borel subset of } \mS^{n-1}. $$
Here ${\cal{H}}^{n-1}(\cdot)$ stands for the $(n-1)$-dimensional Hausdorff measure.

Let $V(K_1,\dots,K_n)$ denote the $n$-dimensional mixed volume of
$n$ convex bodies $K_1,\dots, K_n$ in $\R^n$. We write $V(K_1[m_1],\dots, K_r[m_r])$ for the
mixed volume of the bodies $K_1,\dots, K_r$ where each $K_i$ is repeated $m_i$ times and
$m_1+\dots+m_r=n$. In particular, $V(K[n])=V_{n}(K)$, the $n$-dimensional Euclidean volume of $K$.  

Let $S(K_1,\dots, K_{n-1},\cdot)$ be the {\it mixed area measure} for bodies $K_1,\dots, K_{n-1}$ defined by
$$V(L,K_1,\dots, K_{n-1})=\frac{1}{n}\int_{\mS^{n-1}}h_LdS(K_1,\dots, K_{n-1},\cdot)$$
for any compact convex set $L$.
In particular, when the $K_i$ are polytopes the mixed area measure $S(K_1,\dots, K_{n-1},\cdot)$ has finite support and
for every $u\in\mS^{n-1}$ we have
\begin{equation}\label{e:mixed-area}
S(K_1,\dots, K_{n-1},u)=V(K_1^u,\dots, K_{n-1}^u),
\end{equation}
where $V(K_1^u,\dots, K_{n-1}^u)$ is the $(n-1)$-dimensional mixed volume of the faces $K_i^u$ translated 
the the subspace orthogonal to $u$,
see \cite[Sec 5.1]{Sch1}.

Finally, for $u\in\mS^{n-1}$ the orthogonal projection of a set $A\subset\R^n$ onto the subspace  $u^{\bot}$
orthogonal to $u$ is denoted by $A|u^{\bot}$.


\section{Proof of \rt{simplex0}}

In this section we give a proof of \rt{simplex0}. As mentioned in the introduction, it is enough to prove it for $r=2$ in which case we write the Bezout inequality as
\begin{equation}\label{e:BI-2}
V(L,M,K[n-2])V_n(K)\leq V(L,K[n-1])V(M,K[n-1]).
\end{equation}
We assume that $L$, $M$ are arbitrary convex bodies and $K$ is a polytope in $\R^n$.

We need to set up additional notation. Let $K$ be defined by inequalities
$$K=\bigcap_{j=1}^N\{x\in\R^n : \la x,u_j\ra\leq h_K(u_j)\},$$
where $u_j$ are the outer normals to the facets of $K$ (in some fixed order) and $N$ is the number of facets of $K$.
Denote by $K_{t,i}$ the polytope obtained by moving the $i$-th facet of $K$ by $t$, that is
$$K_{t,i}=\bigcap_{j=1\atop j\neq i}^{N}\{x\in\R^n : \la x,u_j\ra\leq h_K(u_j)\}\bigcap \{x\in\R^n : \la x,u_i\ra\leq h_K(u_i)+t\}.$$
By abuse of notation we let $K_t$ denote $K_{t,N}$.

\begin{Lemma}\label{L:support} Let $K$ and $K_t$ be as above. Then there exists
$\delta=\delta(K)$ such that the following supports are equal
$$\supp S(K_t[r],K[n-1-r],\cdot)=\supp S(K,\cdot)$$
for any $0\leq r\leq n-1$ and any $t\in(-\delta, \delta)$.
\end{Lemma}

\begin{pf}
By \re{mixed-area} it is enough to show that $V(K_t^u[r],K^u[n-1-r])=0$ if and only if $V_{n-1}(K^u)=0$,
that is $K^u$ is not a facet of $K$. Indeed, by choosing $\delta$ small enough we can ensure that
$K_t$ has the same facet normals as $K$ and so $\dim K_t^u=n-1$ whenever $K^u$ is a facet of $K$.
In this case $V(K_t^u[r],K^u[n-1-r])>0$. 

Conversely, assume $K^u$ is a face of $K$ of dimension less than $n-1$. As before, for small enough
$t$ the face $K_t^u$ also has dimension less than $n-1$. First, suppose $K^u$
is not contained in the moving facet $F=K\cap H_K(u_N)$. Then $h_K(u)=h_{K_t}(u)$ and
so $K^u\subseteq K_t^u$ for $t\geq 0$ and $K^u\supseteq K_t^u$ for $t<0$. Then, by the monotonicity 
of the mixed volume, if $t\geq 0$ then
$$0\leq V(K_t^u[r],K^u[n-1-r])\leq V_{n-1}(K_t^u)=0,$$
and so $V(K_t^u[r],K^u[n-1-r])=0$. The case $t<0$ is similar.

Now suppose $K^u$ is contained in the moving facet $F$. Then $K^u\subseteq H_K(u)\cap H_K(u_N)$
and $K_t^u\subseteq H_{K_t}(u)\cap H_{K_t}(u_N)$. This shows that $K^u$ and $K_t^u$ are 
contained in two affine $(n-2)$-dimensional subspaces which are translates of the same linear subspace
of dimension $n-2$. Therefore, for any collection of line segments $(L_1,\dots, L_{n-1})$,
where $L_i\subset K_t^u$ for $1\leq i\leq r$ and $L_i\subset K^u$  for $r+1\leq i\leq n-1$, the $L_i$ have linearly
dependent directions. The latter implies that  $V(K_t^u[r],K^u[n-1-r])=0$ by \cite[Theorem 5.1.7]{Sch1}.

\end{pf}

\begin{Prop} \label{P:r=1}
Let $K, P$ be $n$-polytopes with the following properties:
\begin{enumerate}
 \item $\supp S(P,\cdot)=\supp S(K,\cdot)$,
 \item there exists a constant $\lambda>0$ such that $V(L,P[n-1])\leq \lambda V(L,K[n-1])$ for all convex bodies $L$,
 \item $V(K,P[n-1])=\lambda V_n(K)$.
 \end{enumerate}
Then,
$$S(P,\cdot)=\lambda S(K,\cdot)\ .$$
\end{Prop}

\begin{pf} As before, let $\{u_1,\dots, u_N\}$ be the outer normals to the facets of $K$. 
By assumption (1) they are the outer normals to the facets of $P$ as well.
Fix $1\leq i\leq N$ and let $L=K_{s,i}$ be the polytope obtained from $K$ by moving its $i$-th facet by a small 
number $s\in(-\delta_i,\delta_i)$ as in \rl{support}. 
 
By assumption (2), for any $s\in(-\delta_i,\delta_i)$ we have
$$V(K_{s,i},P[n-1])\leq \lambda V(K_{s,i},K[n-1]).$$
Consider the function
$$F(s)=\lambda V(K_{s,i},K[n-1])-V(K_{s,i},P[n-1]).$$
Then $F(s)\geq 0$ and $F(0)=0$. Below we show that $F(s)$ is, in fact, linear on $(-\delta_i,\delta_i)$. But then $F(s)$ is identically zero on $(-\delta_i,\delta_i)$, which implies that 
\begin{equation}\label{e:equality}
V(K_{s,i},P[n-1])=\lambda V(K_{s,i},K[n-1])
\end{equation}
 for all $s\in(-\delta_i,\delta_i)$.
We claim that this also implies that
\begin{equation}\label{e:claim2}
S(P,u_i)=\lambda S(K,u_i),
\end{equation}
and since $i$ is chosen arbitrarily and the supports of the two measures
are equal, the statement of the proposition follows.
 
Now we show that $F(s)$ is linear and then prove that \re{equality} implies \re{claim2}.
Since the polytopes $P$ and $K$ have the same set of facet normals
$\{u_1,\dots,u_N\}$, we obtain:
\begin{eqnarray}\label{e:1}
nV(K_{s,i}, P[n-1])&=&\sum_{j=1}^Nh_{K_{s,i}}(u_j)V_{n-1}(P^{u_j})\nonumber\\
&=&\sum_{j=1}^{N}h_{K}(u_j)V_{n-1}(P^{u_j})+(h_{K}(u_i)+s)V_{n-1}(P^{u_i})\nonumber\\
&=&nV(K,P[n-1])+sV_{n-1}(P^{u_i})\nonumber\\
&=&n\lambda V_n(K)+sV_{n-1}(P^{u_i}).
\end{eqnarray}
Similarly, 
\begin{equation}\label{e:2}
nV(K_{s,i},K[n-1])=nV_n(K)+sV_{n-1}(K^{u_i}).
\end{equation}
Substituting \re{1} and \re{2} into the definition of $F(s)$ and using assumption (3), we see that
$F(s)=\lambda s$ for some $\lambda$, that is  $F(s)$ is linear.

It remains to show that \re{equality} implies \re{claim2}. Since $F(s)$ is identically zero we have $\lambda=0$,
which translates to
$$V_{n-1}(P^{u_i})=\lambda V_{n-1}(K^{u_i}).$$
But that is precisely what \re{claim2} is stating, which completes the proof of the proposition.

\end{pf}

\begin{Lemma}\label{L:any r}
Let $K$  be an $n$-polytope satisfying \re{BI-2} for all bodies $L$ and for all $M=K_t$ 
where $t\in(-\delta, \delta)$ as in \rl{support}. Then 
$$S(K_t[r],K[n-1-r],\cdot)=\frac{V(K_t,K[n-1])^r}{V_n(K)^r}S(K,\cdot)$$
for all $0\leq r\leq n-1$ and all $t\in(-\delta, \delta)$.
\end{Lemma}

\begin{pf} 
For $0\leq r \leq n-1$, set $P_r$ to be the polytope whose surface area measure equals $S(K_t[r],K[n-r-1],\cdot)$ and let $\lambda:=V(K_t,K[n-1])/V_n(K)$. For each $r$ the existence and uniqueness of $P_r$ is ensured by the Minkowski Existence and Uniqueness Theorem (see \cite[Sections 7.1, 7.2]{Sch1}). We need to prove that 
\begin{equation}
S(P_r,\cdot)=\lambda^rS(K,\cdot),\qquad r=0,1,\dots,n-1.\label{e:proof} 
\end{equation}
Note that by \rl{support}, we have:
\begin{equation}
\supp S(P_r,\cdot)=\supp S(K,\cdot),\qquad r=1,\dots,n-1.\label{e:c1} 
\end{equation}
We prove \re{proof} by induction on $r$. The case $r=0$ is trivial. For the case $r=1$ we apply 
\rp{r=1} with $P=P_1$. Indeed, by our assumption, \re{BI-2} is satisfied for $M=K_t$ and becomes equality when $L=K$.
Thus the conditions (1)--(3) of \rp{r=1} hold and so $S(P_1,\cdot)=\lambda S(K,\cdot)$, as required.

Now assume  \re{proof} holds for $1\leq m\leq r-1$. This is equivalent to the following: 
\begin{equation}\label{e:equivalent}
V(L,P_m[n-1])=\lambda^mV(L,K[n-1]),
\end{equation}
for all convex bodies $L$ and $1\leq m\leq r-1$.
Next fix a convex body  $L\subset \R^n$ and apply the Aleksandrov-Fenchel inequality
\begin{eqnarray*}
&&V(L,P_{r-1}[n-1])^2=V(L,K_t[r-1],K[n-r])^2\\
&=&V(K,K_t,K_t[r-2],K[n-r-1],L)^2\\
&\geq& V(K,K,K_t[r-2],K[n-r-1],L)V(K_t,K_t,K_t[r-2],K[n-r-1],L)\\
&=&V(L,K_t[r-2],K[n-r+1])V(L,K_t[r],K[n-r-1])\\
&=&V(L,P_{r-2}[n-1])V(L,P_r[n-1])\ ,
\end{eqnarray*}
which, by \re{equivalent} with $m=r-2$ and $m=r-1$, gives
$$\lambda^{2(r-1)}V(L,K[n-1])^2\geq \lambda^{r-2}V(L,K[n-1])V(L,P_r[n-1]).$$
Thus
\begin{equation}\label{e:c2}
 V(L,P_r[n-1])\leq \lambda^rV(K,P_r[n-1]).
\end{equation}
Furthermore, using \re{equivalent} for $m=r-1$, we get:
\begin{eqnarray}
V(K,P_r[n-1])&=&V(K,K_t[r],K[n-1-r])\nonumber\\
&=&V(K_t,K_t[r-1],K[n-r])\nonumber\\
&=&V(K_t,P_{r-1}[n-1])\nonumber\\
&=&\lambda^{r-1}V(K_t,K[n-1])\nonumber\\
&=&\frac{V(K_t,K[n-1])^{r-1}}{V_n(K)^{r-1}}V(K_t,K[n-1])=\lambda^rV_n(K).\label{e:c3}
\end{eqnarray}
Now, as in the case of $r=1$, \re{c1}, \re{c2}, \re{c3} together with \rp{r=1}, show that $S(P_r,\cdot)=\lambda^rS(K,\cdot)$, which completes the proof of the lemma.
\end{pf}

Now we are ready to prove the main theorem which implies \rt{simplex0}.

\begin{Th}\label{T:simplex}
Let $K$ be an $n$-polytope in $\R^n$. Suppose that
\begin{equation}\label{e:BI-3}
V(L,M,K[n-2])V_n(K)\leq V(L,K[n-1])V(M,K[n-1])
\end{equation}
holds for all convex bodies $L$ and $M$ in $\R^n$. Then $K$ is a simplex.
\end{Th}

\begin{pf} Let $K_t$ be the polytope obtained by moving one of the facets of $K$ for $t$ small enough.
Then \rl{any r} with $r=n-1$ implies that the surface area measures of $K_t$ and $K$ are proportional, and
hence, $K_t$ is homothetic to $K$. 

We may assume that one of the vertices of $K$ not lying on the moving facet
is at the origin, so $K_t=\lambda K$ for some $\lambda\neq 1$. For every vertex $v$ in $K$, $\lambda v$ must be
a vertex of $\lambda K$. Therefore, the origin is the only vertex of $K$ not lying on the moving facet. In other
words, $K$ is the cone over the moving facet. But since the facet was chosen arbitrarily, for every vertex $v$ the polytope $K$ is the convex hull of $v$ and the facet not containing $v$. This implies that $K$ is a simplex.
\end{pf}


\section{Proof of \rt{last0}}

Recall that a boundary point $y\in\bd K$ is \emph{strict} 
if it does not belong to any segment contained in $\bd K$.  Note that points with positive Gaussian curvature
and, more generally, regular exposed points are strict points (see \cite{Sch1} for the definitions).
Clearly the boundary of a polytope does not  contain any strict points,  but there are other convex bodies having this property (for example, a cylinder). 

As before it is enough to prove \rt{last0} in the case of $r=2$. It follows from the theorem below.

\begin{Th}\label{T:last}
Let $K$ be a convex body whose boundary contains at least one strict point. 
Then there exist convex bodies $L$ and $M$ such that
\begin{equation}\label{e:notBI}
V(L,M,K[n-2])V_n(K)> V(L,K[n-1])V(M,K[n-1]).
\end{equation}
\end{Th}
\begin{pf}
First let us fix some notation. For $a>0$ and $u\in\mS^{n-1}$, define the closed half-spaces:
$$H_a^+(u)=\{x\in\mathbb{R}^n:\langle x, u\rangle \geq a\}\quad \textnormal{and}\quad H_a^-(u)=\{x\in\mathbb{R}^n:\langle x, u\rangle \leq a\}.$$
Also set $H_a(u):=H^+_a(u)\cap H_a^-(u)$. With this notation, the supporting hyperplane of $K$ whose unit normal vector is $u$, can be written as $H_{h_K(u)}(u)$. 

Let $y$ be a strict point of $\bd K$  and $u$ be a normal vector of $K$ at $y$.  
Choose $v\in \mS^{n-1}$, such that $y|v^{\bot}\in \relint(K|v^{\bot})$, where $\relint(K|v^{\bot})$ denotes the relative interior of the body $K|v^{\bot}$ in $v^\bot$.
We claim that there exists $\varepsilon>0$, such that 
\begin{equation}\label{e:last1}
\Big(K\cap H^-_{h_K(u)-\varepsilon}(u)\Big)|v^{\bot}=K|v^{\bot}. 
\end{equation}
To see this, assume that \re{last1} is not true for all $\varepsilon>0$. This means that for any $\varepsilon>0$, there exists a point $x_{\varepsilon}\in \bd K$, such that $x_{\varepsilon}|v^{\bot}\in\bd(K|v^{\bot})$ and $x_{\varepsilon}\in H^+_{h_K(u)-\varepsilon}(u)$. Let $x_0$ be an accumulation point of the set $\{x_{\varepsilon}:\varepsilon>0\}$. Then, by compactness, $x_0\in\bd K$, $x_0|v^{\bot}\in\bd(K|v^{\bot})$, and $x_0\in H_{h_K(u)}(u)$ (because $x_0\in H^+_{h_K(u)}(u)$ and $x_0\in K$).
Note that, since $x_0|v^{\bot}\in \bd(K|v^{\bot})$ and $y|v^{\bot}\in \relint(K|v^{\bot})$, we have $x_0\neq y$. It follows that the segment $[x_0,y]$ is contained in a supporting hyperplane of $K$, thus $[x_0,y]\subseteq\bd K$, which contradicts the assumption that $y$ is strict. Hence, \re{last1} holds for some $\varepsilon>0$.

Next, set $K_{\varepsilon}:=K\cap H^-_{h_K(u)-\varepsilon}(u)$. Clearly, $h_{K_{\varepsilon}}\leq h_K$. We claim that there exists an open subset $\beta\subset \bd K\setminus\bd K_{\varepsilon}$, such that  $y\in \beta$ and 
\begin{equation}\label{e:last2}
h_{K_{\varepsilon}}(u)<h_K(u),\qquad \textnormal{for all }u\in \sigma(K,\beta). 
\end{equation}
Suppose not. Then for any $\delta$-neighborhood 
$\beta_\delta=( \bd K\setminus\bd K_{\varepsilon})\cap B(y,\delta)$ of $y$
there exists a unit vector $u_\delta\in\sigma(K,\beta_\delta)$ such that $h_K(u_\delta)=h_{K_\varepsilon}(u_\delta)$.
In other words, there exist
points $y_{\delta}\in\beta_\delta$ and $x_{\delta}\in\bd K_{\varepsilon}$ lying in the same hyperplane $H_K(u_\delta)$.
But then, by compactness,
there exist a point $x\in\bd K_{\varepsilon}$ and a unit vector $u$, which is normal for $K$ at $y$ and at $x$. 
This shows again that the points $y$ and $x$ of $K$ lie in the same supporting hyperplane $H_K(u)$, 
thus $[y,x]$ is a boundary segment of $K$, which contradicts our assumption. Therefore, \re{last2} holds for some open set $\beta\subseteq\bd K\setminus\bd K_{\varepsilon}$.

Note, furthermore, that $\tau(K,\sigma(K,\beta))\supseteq \beta$, thus ${\cal{H}}^{n-1}\big(\tau(K,\sigma(K,\beta))\big)>0$, which shows that 
\begin{equation}\label{e:last3}
S(K,\sigma(K,\beta))>0. 
\end{equation}
Now we are ready to exhibit examples of compact convex sets $L$ and $M$ satisfying \re{notBI}.
Set $L=[-v,v]$ and $M=K_{\varepsilon}$. Then, by (5.3.23) in \cite[p. 294]{Sch1} and applying \re{last1} we obtain
$$V(L,M,K[n-2])=V(K_{\varepsilon}|v^{\bot},K|v^{\bot}[n-2])=V_{n-1}(K|v^{\bot})=V(L,K[n-1]).$$
 
On the other hand, by \re{last2} and \re{last3}, we have:
\begin{eqnarray*}V(M,K[n-1])=V(K_{\varepsilon},K[n-1])&=&\frac{1}{n}\int_{S^{n-1}}h_{K_{\varepsilon}}dS(K,\cdot)\\ 
&<&\frac{1}{n}\int_{S^{n-1}}h_{K}dS(K,\cdot)=V_n(K).
\end{eqnarray*}
This shows that $$V(L,M,K[n-2])V_n(K)>V(L,K[n-1])V(M,K[n-1]),$$
as asserted.
\end{pf}
\begin{Rem}
One might ask the following: 
If $K$ is a convex body whose boundary contains at least one strict point $x$, is it true that $\bd K$ has an open neighborhood that does not contain any line segments, i.e. $K$ is strictly convex in a neighborhood of $x$? If yes, this would simplify the proof of \rt{last} considerably. 
The following simple 3-dimensional example shows, however, that this is not the case. Take $K$ equal to
$$\{x\in\mathbb{R}^3:x_3\leq 1\}\bigcap \textnormal{conv}\Big(\{(0,x_2,x_3)\in\mathbb{R}^3:x_3=x_2^2\}\cup\{(x_1,0,x_3)\in\mathbb{R}^3:x_3=x_1^2\}\Big).$$
Then the origin is a strict point of the boundary of $K$, but no neighborhood of the origin is strictly convex.
\end{Rem}

\end{document}